\documentclass{article}
\usepackage{amsmath,amsfonts,mathrsfs,amssymb,amscd,comment,enumerate,amsxtra,url,tikz-cd,tcolorbox,graphicx, setspace, mathtools,stmaryrd, mathabx,physics,hyperref}
\usepackage{amsthm}

\usepackage{cleveref} 
\hypersetup{colorlinks=true, linkcolor=blue, citecolor=blue, urlcolor=blue}
\input{xypic}
\xyoption{all}
\theoremstyle{plain}
\newtheorem{theorem}{Theorem}[section]
\newtheorem{proposition}[theorem]{Proposition}
\newtheorem{lemma}[theorem]{Lemma}
\newtheorem{corollary}[theorem]{Corollary}

\newtheorem{introthm}{Theorem}

\theoremstyle{definition}
\newtheorem{definition}[theorem]{Definition}
\newtheorem{remark}[theorem]{Remark}

\crefname{theorem}{Theorem}{Theorems}
\crefname{proposition}{Proposition}{Propositions}
\crefname{lemma}{Lemma}{Lemmas}
\crefname{corollary}{Corollary}{Corollaries}
\crefname{definition}{Definition}{Definitions}
\crefname{remark}{Remark}{Remarks}
\crefname{example}{Example}{Examples}
\crefname{equation}{}{}
\Crefname{theorem}{Theorem}{Theorems}
\Crefname{proposition}{Proposition}{Propositions}
\Crefname{lemma}{Lemma}{Lemmas}
\Crefname{corollary}{Corollary}{Corollaries}
\Crefname{definition}{Definition}{Definitions}
\Crefname{remark}{Remark}{Remarks}
\Crefname{example}{Example}{Examples}

\crefname{introthm}{Theorem}{Theorems}

\newcommand{\Pic}{\textnormal{Pic}}
\newcommand{\Ext}{\textnormal{Ext}}

\newcommand{\Hom}{\operatorname{Hom}}

\newcommand{\Quot}{\textnormal{Quot}}
\newcommand{\tf}{\textnormal{tf}}

\newcommand{\rk}{\textnormal{rank}}

\providecommand{\Bun}{\operatorname{Bun}}
\providecommand{\Coh}{\operatorname{Coh}}

\providecommand{\End}{\operatorname{End}}

\providecommand{\homdim}{\operatorname{hom}}
\providecommand{\extdim}{\operatorname{ext}}

\providecommand{\HHom}{\mathcal{H}\!om}

\providecommand{\tors}{\operatorname{tors}}
\title{Counting maximal subbundles on singular curves}
\author{Arkamouli Debnath}
\date{\today}
\begin{document}
\maketitle
\begin{abstract}
    Counting maximal degree subbundles of a given general stable vector bundle $E$ on a smooth projective curve $C$ has been a classical enumerative question. From Segre's work in 1889 for the rank $2$ case, all the way to Holla's work in 2004 for the general case, the question is settled for smooth curves, with an explicit formula for the count. In this paper, we investigate the same question for singular curves. The presence of torsion sheaves at the singular points, and torsion free but not locally free sheaves, makes the boundary of the relevant Quot scheme more complicated in our case. We show that for a general stable bundle on a curve with one ADE singularity, in every rank, and in rank $2$ on a curve with locally planar singularities, the relevant Quot scheme parametrizes only vector bundle quotients, is finite and reduced, and the count matches Holla's number.
\end{abstract}

\tableofcontents
\section{Introduction}\label{sec:introduction}

Let $C$ be an integral projective complex curve of arithmetic genus
$g\ge2$, let $E$ be a vector bundle of rank $r$ and degree $d$ on $C$,
and fix an integer $0<k<r$.

\begin{definition}\label{def:maximal-subbundle}
A \emph{subbundle} of $E$ is a sub-vector bundle $F\subset E$ whose quotient
$E/F$ is also locally free. Let
$e_{\max}(E)$ be the maximal degree of a rank-$k$ subsheaf of $E$. A
\emph{maximal subbundle} of rank $k$ is a subbundle of rank $k$ and maximal
degree $e_{\max}(E)$.
\end{definition}

\begin{tcolorbox}
\textbf{Question.} Let $E$ be a general stable bundle of rank $r$ and
degree $d$ on $C$.
\begin{enumerate}
\item[(i)] For which $r,d,k$ does $E$ have only finitely many maximal
subbundles of rank $k$?
\item[(ii)] When it does, is there a closed formula for their number?
\end{enumerate}
\end{tcolorbox}

\subsection*{Smooth curves}
For smooth curves $C$, this is a very classical question. When the rank $r=2$, the answer goes back to Segre's work on
ruled surfaces \cite{Seg89}: a general stable bundle of rank two and
degree $d\equiv g-1\pmod 2$ has exactly $2^g$ maximal line subbundles.
An interesting observation is that, this number depends only on the genus, and not on the curve or on the
general bundle chosen.

In arbitrary rank, the rank-$k$ subsheaves of $E$ of degree $e$ are
the points of the Quot scheme $\Quot_{d-e,r-k}(E)$ of quotients of $E$
of rank $r-k$ and degree $d-e$. At a point $0\to F\to E\to Q\to0$ its
tangent and obstruction spaces are $\Hom(F,Q)$ and $\Ext^1(F,Q)$, so
on a smooth curve each of its components has dimension at least
\begin{equation}\label{eq:intro-expected-dimension}
 D_e=\chi(F^\vee\otimes Q)=kd-re-k(r-k)(g-1).
\end{equation}
Hirschowitz proved that, for a general stable bundle, the maximal
degree of a rank-$k$ subbundle is the largest $e$ with $D_e\ge0$
\cite{Hir88}. A finite answer can therefore be expected only in the
\emph{critical case}
\begin{equation}\label{eq:intro-critical}
 kd\equiv k(r-k)(g-1)\pmod r,\qquad
 e_*=\frac{kd-k(r-k)(g-1)}{r},
\end{equation}
in which case $D_{e_*}=0$. In this case Holla proved that a general stable
bundle has finitely many maximal subbundles, and computed their number
$m(r,d,k,g)$ by means of Gromov--Witten invariants of Grassmannians
and the formula of Vafa and Intriligator
\cite[Proposition~4.1 and Theorem~4.2]{Hol04}; his formula is recalled
in \eqref{eq:ade-holla-count}. For $k=1$ it gives $r^g$
\cite[Corollary~4.3]{Hol04}, which is Segre's number when $r=2$. Like
Segre's number, $m(r,d,k,g)$ depends only on the numerical data. We
refer to the introduction of \cite{Hol04} for the history of the
problem.

\subsection*{Singular curves}
Smooth curves degenerate to singular ones, and understanding how a
phenomenon on smooth curves remains or changes at the boundary is a
central theme in algebraic geometry. Since $m(r,d,k,g)$ does not
depend on the smooth curve, it is natural to ask whether a general
stable bundle on an integral singular curve of arithmetic genus $g$
has finitely many maximal subbundles, and if so, does that match Holla's count.

This is however not a formal consequence of the smooth case. On a singular
curve, a saturated subsheaf of $E$ has a torsion-free quotient which
need not be locally free, so a subsheaf of maximal degree need not be
a subbundle. For the Question to have the same answer as on a smooth
curve, one must show that for a general $E$ no such subsheaf occurs at
the critical degree, and that the Quot scheme is finite and reduced
there; only then is its number of points a number of subbundles, and then it remains to identify this number.

For irreducible nodal curves, Bhosle proved that the tensor product of
two general stable bundles is nonspecial, that the maximal degree is
given by Hirschowitz's formula, and that in the critical case the
maximal subsheaves of a general stable bundle are finitely many
subbundles, forming a smooth Quot scheme
\cite[Theorems~1 and~2, Propositions~2.8 and~3.7(3)]{Bho06}.

To our knowledge, the number of maximal subbundles in the critical
case has not previously been determined for a general stable bundle
on a fixed integral singular curve of arithmetic genus at least two.
Teixidor computes the smooth-curve numbers for $k=1$ and for
$(r,k)=(4,2)$ by degenerating to a reducible curve with one node
\cite{Tei04}, and Ballico bounds the maximal degree of subsheaves on
integral curves \cite{Bal01}.
 
\subsection*{Results}
We answer the Question in the critical case for two classes of
singular curves. These are curves with one singular point
of ADE type (\cref{def:ade-singularity}), in all ranks, and curves with locally planar singularities, for the rank $2$ case.

\begin{introthm}\label{thm:intro-quot}
Let $C$ be an integral projective complex curve of arithmetic genus
$g\ge2$ with exactly one singular point, of ADE type, and assume the
critical congruence \eqref{eq:intro-critical}. There is a nonempty
open substack of the stack of vector bundles of rank $r$ and degree
$d$ on $C$, consisting of stable bundles $E$, such that $e_*$ is the
maximal degree of a rank-$k$ subsheaf of $E$ and the Quot scheme
$\Quot_{d-e_*,r-k}(E)$ is nonempty, finite and reduced. Each of its
points is a sequence $0\to F\to E\to Q\to0$ of stable vector bundles
with
\[
 \Hom_C(F,Q)=\Ext^1_C(F,Q)=0 .
\]
\end{introthm}

\begin{introthm}\label{thm:intro-count}
In the situation of \cref{thm:intro-quot}, a general stable bundle of
rank $r$ and degree $d$ on $C$ has exactly $m(r,d,k,g)$ maximal
subbundles of rank $k$. In particular, it has $r^g$ maximal line
subbundles when $k=1$.
\end{introthm}

In rank two the ADE hypothesis can be replaced by local planarity.

\begin{introthm}\label{thm:intro-planar}
Let $C$ be an integral projective complex curve of arithmetic genus
$g\ge2$ with a locally planar singularity, and let
$d\equiv g-1\pmod 2$. A general stable bundle $E$ of rank two and
degree $d$ on $C$ has exactly $2^g$ maximal line subbundles. They are
the points of the Quot scheme $\Quot_{d-e_*,1}(E)$, where
$e_*=(d-g+1)/2$, which is finite and reduced.
\end{introthm}

\Cref{thm:intro-quot,thm:intro-count} are
\cref{thm:ade-critical-quot,thm:ade-counting} in
\cref{sec:ade-critical}, and \cref{thm:intro-planar} is
\cref{thm:rank-two-planar} in \cref{subsec:rank-two-planar}. Thus, on
these curves, a subsheaf of maximal degree of a general stable bundle
is a subbundle in the critical case, and the number of such subbundles
is the one found by Segre and Holla on smooth curves.

For curves with one ADE singularity we also obtain Hirschowitz's
formula for the maximal degree (\cref{cor:ade-hirschowitz}), the
nonspeciality of the tensor product of two general stable bundles
(\cref{prop:ade-nonspeciality}), and a restriction on the local
structure of the quotient of an arbitrary saturated subsheaf of a
general bundle (\cref{prop:ade-generic-geometry}), which is made explicit
for $A_n$ singularities (\cref{cor:an-explicit-penalty}).

The proofs are outlined in \cref{subsec:counting-background}.

\paragraph{Acknowledgments. } The author thanks Parvez Rasul for introducing him to Quot schemes on singular curves. The author thanks Michele Pernice for useful discussions. 

\paragraph{AI Statement. } The author used Claude Fable 5.1 and ChatGPT 6.1 for literature search on existing results. AI was also used to make the explicit computations in the appendix for "local penalties" for the $A_n$ singular curves. The author has independently verified the arguments and all mathematical claims in the paper and takes full responsibility for their correctness.

\section{Strategy of the proof}\label{subsec:counting-background}

Throughout, we put $s=r-k$, we assume the critical congruence
\eqref{eq:intro-critical}, and we call $\Quot_{d-e_*,s}(E)$ the
\emph{critical Quot scheme} of $E$; it parametrizes all quotients of
$E$ of rank $s$ and degree $d-e_*$, with or without torsion. The
proofs of \cref{thm:intro-quot,thm:intro-count,thm:intro-planar} rest
on two general results, in which no particular singularity appears.

The first is a degeneration argument (\cref{prop:count-transport}).
Let $C$ be integral and locally planar, and let $E$ be a vector bundle
whose critical Quot scheme is nonempty and finite, each of its points
being a sequence of vector bundles $0\to F\to E\to Q\to0$ with
$\Hom_C(F,Q)=\Ext^1_C(F,Q)=0$. Then this scheme is reduced and has
exactly $m(r,d,k,g)$ points: over a smoothing of $C$ it deforms to a
finite \'etale cover, whose degree is computed on a smooth fibre by
Holla's theorem.

The second is a criterion for a general bundle to satisfy these
hypotheses, valid on every integral projective curve
(\cref{prop:boundary-generic-geometry}). The map sending a surjection
$E\twoheadrightarrow Q$ onto a torsion-free sheaf to $Q$ is smooth
(\cref{lem:quotient-smooth}), so the dimension of a family of
subsheaves of varying bundles $E$ can be read off from the dimension
of the family of their quotients, and compared with that of the stack
of all bundles. The outcome is that two conditions suffice:
\begin{enumerate}
\item[(B)] \emph{Boundary estimate.} Whenever $D_e\le0$, the
torsion-free sheaves of rank $s$ and degree $d-e$ that are not
locally free form a locus of dimension at most $s^2(g-1)-1$, one
less than the dimension of the stack of vector bundles of rank $s$.
\item[(N)] \emph{Nonspeciality.} A general pair $(F,Q)$ of vector
bundles of ranks $k,s$ and degrees $e_*,d-e_*$ satisfies
$\Hom_C(F,Q)=0$.
\end{enumerate}
If (B) and (N) hold, then for a general bundle $E$ the degree $e_*$ is
the maximal degree of a rank-$k$ subsheaf, and the critical Quot
scheme is nonempty, finite and reduced, with points as above. On a
locally planar curve the first result then gives the count. It
remains to verify (B) and (N), which we do in two situations.

\medskip\noindent
\emph{Locally planar curves in rank two}
(\cref{subsec:rank-two-planar}). Here $k=s=1$, so both conditions
concern rank-one sheaves and follow from known rank-one geometry:
(B) from the density of line bundles in the compactified Jacobian
\cite[Theorem~A(ii)]{MRV15} (\cref{lem:planar-rank-one-boundary}), and
(N) from the fact that a general line bundle of degree $g-1$ has no
sections (\cref{lem:planar-line-nonspeciality}).

\medskip\noindent
\emph{One ADE singularity}
(\cref{sec:local-types,sec:ade-critical}). Here $s$ is arbitrary,
and no global moduli space of the boundary is used. Instead, we
stratify non-locally-free sheaves by their \emph{local type}, the
isomorphism class $M$ of the completed stalk at the singular point.
Each type has codimension at least an explicit
\emph{penalty} $c_R(M)$ (\cref{prop:ade-type-dimension}). The
penalty is positive for nonfree $M$ because $R$ is Gorenstein
(\cref{prop:local-gorenstein-penalty}), and there are finitely many
types in each rank \cite{GK85}, which gives (B)
(\cref{cor:ade-boundary}). Condition (N) is no longer straightforward as in the rank $2$ case; it
follows from the argument of Hirschowitz and Bhosle
\cite{Hir88,Bho06}, whose geometric inputs are consequences of (B)
(\cref{prop:ade-nonspeciality}).

\medskip
The count uses only the positivity of the penalty.
\Cref{app:an-example} computes it explicitly for $A_n$ singularities,
which makes it a numerical restriction on the non-locally-free
quotients of a general bundle when $D_e>0$
(\cref{prop:ade-generic-geometry,cor:an-explicit-penalty}).
\vspace{2mm}
\begin{remark}[Comparison with \cite{Reg82}]\label{rem:rego}
For an integral curve lying on a smooth surface, that is, with
locally planar singularities \cite[\S1.6]{MRV15}, Rego proved that
every torsion-free sheaf is a limit of vector bundles
\cite[(0.1) and Theorem~2.5]{Reg82}. This would give (B) on all such
curves, in every rank, with no analysis of local types. We do not use
this result: our proof of (B) is independent of it, and it gives an
explicit lower bound $c_R(M)$ for the codimension of each local type
rather than density alone, which is the source of the refined bound
in \cref{prop:ade-generic-geometry}.
\end{remark}  

\section{Common setup and the quotient map}\label{sec:quotient-map}

\textbf{Conventions and Notations: }In this section $C$ is an integral projective complex curve, with
$g=p_a(C)$, and $\nu:Y\to C$ is its normalization. Degrees of coherent sheaves are defined by
\[
 \deg B=\chi(B)-\rk(B)\chi(\mathcal O_C),
\]
where $\rk$ is the generic rank. Degree is additive in short exact
sequences; a torsion sheaf has degree equal to its length.
We fix the rank and degree, but not the determinant.

Write $\Bun_a^b(C)$ and $\Coh^{\tf}_{a,b}(C)$ for the stacks
of vector bundles and torsion-free sheaves of rank $a>0$ and degree $b$.
Both are open substacks of the algebraic stack of coherent sheaves. For the torsion-free locus,
openness follows from openness of purity
\cite[Proposition~2.3.1]{HL}, since $C$ is integral.
All dimensions below are stack dimensions, as opposed to corresponding coarse space dimension.
The symbols $\homdim$ and $\extdim$ denote dimensions of global
Hom and Ext groups.

Fix $r,d$, $0<k<r$, and put $s=r-k$. We write
$\Quot_{d-e,s}(E)$ for the full scheme of rank-$s$, degree-$(d-e)$
quotients of $E$. Its saturated open parametrizes sequences
\begin{equation}\label{eq:ade-sequence}
 0\longrightarrow F\longrightarrow E\longrightarrow Q
 \longrightarrow0
\end{equation}
with $Q$ torsion-free, where $F$ has rank $k$ and degree $e$.
Put
\begin{equation}\label{eq:ade-expected-dimension}
 D_e=kd-re-ks(g-1).
\end{equation}
For vector-bundle endpoints this is the number $\chi(F^\vee\otimes Q)$,
the expected Quot dimension. The critical congruence is
\begin{equation}\label{eq:critical-congruence}
 kd\equiv k(r-k)(g-1)\pmod r.
\end{equation}
When it holds, set
\begin{equation}\label{eq:critical-degree}
 e_* =\frac{kd-k(r-k)(g-1)}r;
 \qquad D_{e_*}=0.
\end{equation}
\vspace{3mm}

Various forms of the following result are well known; see
\cite{Alp26},
\cite[Chapter~5, Remark~5.5 and \S7]{New78} and
\cite[Lemma~2.4]{Bho06}. Here we give a proof of the general
statement.

\begin{lemma}\label{lem:ade-bun}
For any integral projective complex curve $D$, the stack
$\Bun_a^b(D)$ is smooth and irreducible of dimension
$a^2(p_a(D)-1)$.
\end{lemma}
\begin{proof}
\emph{Smoothness and dimension.}
For a vector bundle $V$ on $D$, the infinitesimal automorphisms,
first-order deformations and obstructions are $H^0$, $H^1$ and $H^2$
of $\End V$ \cite[Lemma~91.6.1, Tag~08VR]{Stk}. Since $D$ is a curve,
$H^2(D,\End V)=0$, so $\Bun_a^b(D)$ is smooth. Its dimension at $V$ is
therefore
\[
 h^1(\End V)-h^0(\End V)=-\chi(\End V)
 =-\deg(\End V)-a^2\chi(\mathcal O_D)=a^2(p_a(D)-1),
\]
by Riemann--Roch.

\emph{Irreducibility in rank one.}
Let $\nu:Y\to D$ be the normalization, and for $x\in D$ let
$\widetilde{\mathcal O}_{D,x}$ be the integral closure of
$\mathcal O_{D,x}$. By \cite[Exercise~II.6.9(a)]{Har77}, there is an
exact sequence
\[
 0\to K\to\Pic(D)\xrightarrow{\ \nu^*\ }\Pic(Y)\to0,
 \qquad
 K=\prod_{x}\widetilde{\mathcal O}_{D,x}^*/\mathcal O_{D,x}^*,
\]
the product being over the singular points. Pullback preserves
degree, so $\Pic^0(D)$ is an extension of the Jacobian $\Pic^0(Y)$ by
$K$. The Jacobian of the smooth curve $Y$ is connected, and so is
$K$, as a quotient of the unit groups of the finite-dimensional
algebras $\widetilde{\mathcal O}_{D,x}/\mathfrak c_x$, where
$\mathfrak c_x$ is the conductor. Hence $\Pic^0(D)$ is a connected
algebraic group, and therefore irreducible. Each $\Pic^m(D)$ is a
translate of it, and $\Bun_1^m(D)$ is a $\mathbb G_m$-gerbe over
$\Pic^m(D)$.

\emph{Irreducibility in rank $a\ge2$.}
We show that all bundles of rank $a$ and degree $b$ are middle terms
of extensions of a line bundle by a fixed bundle, and that these
extensions form an irreducible family.

Fix an ample line bundle $H$ and an integer $n$ with
$b+an\deg H>0$, and consider extensions
\begin{equation}\label{eq:bun-extension}
 0\to (H^{-n})^{\oplus(a-1)}\to V\to L\to0
\end{equation}
with $L$ a line bundle, necessarily of degree $b+(a-1)n\deg H$. Since
$\deg L>\deg H^{-n}$ and $D$ is integral, $\Hom(L,H^{-n})=0$, so
$\dim\Ext^1(L,H^{-n})$ depends only on the degrees, by Riemann--Roch.
Hence the groups $\Ext^1(L,(H^{-n})^{\oplus(a-1)})$ form a vector
bundle over the Picard variety parametrizing $L$. Its total space
$E_n$ is irreducible by the previous rank-one case, and the middle term of the
universal extension defines a morphism $E_n\to\Bun_a^b(D)$.

Every bundle $V$ lies in the image of $E_n$ once $V\otimes H^n$ is
globally generated. Indeed, an $(a-1)$-tuple of sections of
$V\otimes H^n$ gives $a-1$ vectors in each fibre, and at a given
point these are linearly dependent only for a closed set of tuples of
codimension two. As the point varies over the curve, the tuples that
are dependent somewhere still form a proper closed subset. A general
tuple therefore defines a map
$\mathcal O_D^{\oplus(a-1)}\to V\otimes H^n$ that is injective on
every fibre, so its cokernel is a line bundle, also at singular
points. Twisting by $H^{-n}$ exhibits $V$ as in
\eqref{eq:bun-extension}.

Finally, let $U_1,U_2$ be nonempty opens of $\Bun_a^b(D)$, choose
$V_i\in U_i$, and take $n$ large enough that both $V_i\otimes H^n$ are
globally generated. Then both $V_i$ lie in the image of $E_n$, so the
preimages of $U_1$ and $U_2$ are nonempty opens of the irreducible
$E_n$, and hence they meet. Therefore $U_1\cap U_2\ne\emptyset$.
\end{proof}
\vspace{2mm}
Put $\delta=\ell(\nu_*\mathcal O_Y/\mathcal O_C)=g-g(Y)$.
For a torsion-free sheaf $B$ of rank $a$, set
$B^\nu=\nu^*B/\tors$. The natural injection into its normalization
pushforward and Riemann--Roch give
\begin{equation}\label{eq:ade-normalization-degree}
 \deg B^\nu=\deg B-a\delta+
 \ell(\nu_*B^\nu/B).
\end{equation}
\vspace{1.3mm}
\begin{lemma}[Bounded ambient charts]\label{lem:ade-bounded-charts}
In a finite-type family of rank-$r$, degree-$d$ vector bundles,
the degrees of rank-$k$ subsheaves of geometric fibres are
uniformly bounded above. A lower degree bound therefore leaves
only finitely many degrees. For each degree, the saturated
relative Quot is of finite type.
\end{lemma}
\begin{proof}
Let $\mathcal E$ be the given family on $C\times S$ and choose
a smooth point $p\in C$. For each geometric fibre $E_s$, the
spaces $H^0(E_s(-mp))$ form a descending chain in $H^0(E_s)$
with zero intersection, since a nonzero section cannot vanish
to infinite order at $p$. They are therefore zero for $m$ large.

By upper semicontinuity
\cite[Lemma~36.32.1, Tag~0BDN]{Stk}, the loci where
$H^0(E_s(-mp))=0$ form an increasing open cover of $S$.
Quasi-compactness gives a single $m$ for which this vanishing
holds on every geometric fibre.

For any rank-$k$ subsheaf $F\subset E_s$, the inclusion
$F(-mp)\subset E_s(-mp)$ gives $H^0(F(-mp))=0$.
Since $F$ is locally free of rank $k$ at $p$, twisting by
$-mp$ lowers its Euler characteristic by $km$. Thus
\[
 \deg F-km+k(1-g)
 =\chi(F(-mp))=-h^1(F(-mp))\le0.
\]
Hence $\deg F\le k(m+g-1)$ uniformly, and a lower bound
leaves only finitely many possible degrees.

For fixed $e=\deg F$, the quotient has rank $r-k$ and degree
$d-e$, so its Hilbert polynomial is fixed. The corresponding
relative Quot is of finite presentation over $S$
\cite[Lemma~108.6.4, Tag~0DPD]{Stk}. Its saturated locus is open
by openness of purity \cite[Proposition~2.3.1]{HL}, and is
therefore of finite type since the full relative Quot is
Noetherian.
\end{proof}

\begin{definition}[Saturated relative Quot stack]
\label{def:saturated-quot-stack}
For an integer $e$, let $\mathfrak X_e$ be the stack over
$\Bun_r^d(C)$ defined as follows. An object over a scheme $T$ is
a vector bundle $E_T$ of rank $r$ and degree $d$ on $C\times T$
together with a surjection $E_T\twoheadrightarrow Q_T$ onto a
$T$-flat coherent sheaf whose geometric fibres are torsion-free
of rank $s$ and degree $d-e$, as in \eqref{eq:ade-sequence}.
A morphism is an isomorphism of the ambient bundles carrying
kernel to kernel. We write
\[
 \pi_E:\mathfrak X_e\longrightarrow\Bun_r^d(C),
 \qquad(E\twoheadrightarrow Q)\longmapsto E,
\]
for the forgetful morphism.
\end{definition}

The kernel of an object of $\mathfrak X_e$ is flat over the base,
since both the middle term and the quotient are flat; on each
geometric fibre it is a saturated subsheaf of rank $k$ and
degree $e$. The fibre of $\pi_E$ over a family $E_T$ is the open
locus of the relative Quot scheme of $E_T$ where the quotient is
torsion-free \cite[Theorem~2.2.4 and Proposition~2.3.1]{HL}.
Thus $\pi_E$ is representable, and $\mathfrak X_e$ is an
algebraic stack, locally of finite type.

\begin{lemma}[Smoothness of the quotient map]
\label{lem:quotient-smooth}
The morphism
\[
 \pi_Q:\mathfrak X_e\longrightarrow\Coh^{\tf}_{s,d-e}(C),
 \qquad(E\twoheadrightarrow Q)\longmapsto Q,
\]
is smooth of relative dimension
\begin{equation}\label{eq:rho}
 \rho_e=kd-re+rk(g-1).
\end{equation}
Moreover, we have
\[
 s^2(g-1)+\rho_e=r^2(g-1)+D_e.
\]
\end{lemma}
\begin{proof}
The morphism is locally of finite type by
\cite[Lemma~101.17.8, Tag~06U9]{Stk}; its source and target are
locally Noetherian. We check lifting across a small extension
$A'\twoheadrightarrow A$ of local Artinian complex algebras,
with residue field $L$ and kernel $I$, where
$\mathfrak m_{A'}I=0$. Fix a quotient
$E_A\twoheadrightarrow Q_A$ and an $A'$-flat lift $Q_{A'}$. We seek a dotted arrow making the following diagram
commute:
\[
\begin{tikzcd}[column sep=large, row sep=large]
\operatorname{Spec} A
    \arrow[r, "{[E_A\twoheadrightarrow Q_A]}"]
    \arrow[d, hook, "i"']
&
\mathfrak X_e
    \arrow[d, "\pi_Q"]
\\
\operatorname{Spec} A'
    \arrow[r, "{[Q_{A'}]}"']
    \arrow[ru, dashed,
        "{[E_{A'}\twoheadrightarrow Q_{A'}]}" description]
&
\Coh^{\tf}_{s,d-e}(C).
\end{tikzcd}
\]
Write $F_L$ for the closed-fibre kernel.

On an affine open where $E_A$ is free, lift its generating
sections of $Q_A$ to $Q_{A'}$. The resulting map from a free
rank-$r$ module is surjective by Nakayama. It remains to glue these local lifts. Let
$\{U_i\}$ be an affine open cover of $C$ on which
$E_A$ is free, and write
\[
q_i':E_i'\twoheadrightarrow Q_{A'}|_{U_i}
\]
for the local lifts constructed above.

On an overlap $U_{ij}=U_i\cap U_j$, we first
want an isomorphism
\[
\phi_{ij}:E_i'|_{U_{ij}}
\xrightarrow{\sim}E_j'|_{U_{ij}}
\]
reducing to the original identification of
$E_A$, and satisfying
\[
q_j'\circ\phi_{ij}=q_i'.
\]

An isomorphism lifting the original identification
exists because the bundles are locally free and
the overlap is affine. Its failure to commute with
the quotient maps is a section of
\[
\HHom(E_L,Q_L)\otimes_L I.
\]
Applying $\HHom(E_L,-)$ to the exact sequence
\[
0\to F_L\to E_L\to Q_L\to0
\]
gives a surjection
\[
\End(E_L)\longrightarrow\HHom(E_L,Q_L).
\]
This is also surjective on sections over
affine opens. We can therefore correct the
isomorphism so that it commutes with the
quotient maps.

These corrected isomorphisms need not satisfy
the cocycle condition on triple intersections.
Their failure to do so is measured by
infinitesimal automorphisms of the local
quotient lifts that reduce to the identity.

Such automorphisms are of the form $1+\alpha$,
where $\alpha$ factors through $F_L$.
Consequently, their additive sheaf is
\[
\HHom(E_L,F_L)\otimes_L I
=E_L^\vee\otimes F_L\otimes_L I.
\]

The obstruction to satisfying the cocycle
condition therefore lies in
\[
H^2(C_L,E_L^\vee\otimes F_L)\otimes_L I.
\]
This group vanishes because $C_L$ is a curve.
Hence the local lifts glue to a global quotient
$E_{A'}\twoheadrightarrow Q_{A'}$ extending the original quotient
over $A$, and the infinitesimal lifting criterion proves that
$\pi_Q$ is smooth.

The first-order deformations relative to the
fixed quotient $Q_L$, up to isomorphism, are
parametrized by $H^1(C_L,E_L^\vee\otimes F_L)$,
while the infinitesimal automorphisms are
given by $H^0(C_L,E_L^\vee\otimes F_L)$.
The relative
dimension is therefore
\[
\begin{aligned}
\rho_e
&=-\chi(E_L^\vee\otimes F_L)\\
&=kd-re+rk(g-1),
\end{aligned}
\]
where the last equality follows from
Riemann--Roch.

Finally, using $r=k+s$ and
$D_e=kd-re-ks(g-1)$, we obtain
\[
\begin{aligned}
s^2(g-1)+\rho_e
&=kd-re+(s^2+rk)(g-1)\\
&=r^2(g-1)+D_e.
\end{aligned}
\]
\end{proof}

\begin{corollary}[Quotient-locus dimension bound]\label{cor:penalty-codim}
For a locally closed substack
$\mathfrak T\subset\Coh^{\tf}_{s,d-e}(C)$,
\[
 \dim\pi_Q^{-1}(\mathfrak T)
 \le \dim\mathfrak T+\rho_e
 =r^2(g-1)+D_e-\bigl(s^2(g-1)-\dim\mathfrak T\bigr).
\]
The same bound holds for constructible loci, interpreted on
finite-type smooth charts and their reduced closures.
\end{corollary}
\begin{proof}
The bound follows directly from the relative dimension
formula for the smooth morphism $\pi_Q$
and \eqref{eq:rho}.
\end{proof}
\vspace{2mm}
\begin{lemma}[Smoothness at split sequences]\label{lem:split-smooth}
The morphism $\pi_E:\mathfrak X_e\to\Bun_r^d(C)$ is smooth at every
sequence $0\to F\to E\to Q\to0$ with $\Ext^1_C(F,Q)=0$. In
particular, the image of every open neighbourhood of such a sequence
contains a nonempty open substack of $\Bun_r^d(C)$.
\end{lemma}
\begin{proof}
The fibres of $\pi_E$ are open in relative Quot schemes, and the
obstruction to lifting a quotient along a small extension of the base
lies in $\Ext^1_C(F,Q)$ \cite[Lemma~2.2.6]{HL}. When this group
vanishes, $\pi_E$ is smooth at the sequence by the infinitesimal
criterion. The smooth locus of $\pi_E$ is open, and smooth morphisms
are open \cite[Tag~056G]{Stk}.
\end{proof}
\section{Transport of the critical count}\label{sec:transport}

For a smooth curve of genus $g\ge2$, let $m(r,d,k,g)$ be Holla's
number of maximal rank-$k$ subbundles in the critical case
\eqref{eq:critical-congruence}
\cite[Proposition~4.1 and Theorem~4.2]{Hol04}.
Writing $d=ar-b$, $0\le b<r$, and
$\mu_r=\{\zeta\in\mathbb C:\zeta^r=1\}$, his formula is
\begin{equation}\label{eq:ade-holla-count}
 m(r,d,k,g)=\frac{r^{k(g-1)}}{k!}
 (-1)^{(k-1)(bk-(g-1)k^2)/r}
 \sum_{\substack{(\zeta_1,\ldots,\zeta_k)\in\mu_r^k\\
                  \zeta_i\ne\zeta_j\ (i\ne j)}}
 \frac{(\zeta_1\cdots\zeta_k)^{b-g+1}}
 {\bigl(\prod_{i\ne j}(\zeta_i-\zeta_j)\bigr)^{g-1}}.
\end{equation}
The tuples and the pairs in the product are ordered. The sign
exponent is integral by \eqref{eq:critical-congruence}.
For $k=1$, this gives $m(r,d,1,g)=r^g$
\cite[Corollary~4.3]{Hol04}.

\begin{proposition}[Transport of the critical count]
\label{prop:count-transport}
Let $C$ be integral and locally planar, with $g\ge2$, and assume
\eqref{eq:critical-congruence}. Let $E$ be a rank-$r$, degree-$d$
bundle whose full critical Quot is nonempty and finite. Suppose
every point is a sequence of vector bundles and satisfies
\[
 \Hom_C(F,Q)=\Ext_C^1(F,Q)=0.
\]
Then this Quot is reduced and has length $m(r,d,k,g)$.
There exist a flat projective smoothing of $C$ over
$\mathbb C[[t]]$ and a vector-bundle extension of $E$.
For every such smoothing and extension, the full critical
relative Quot is finite \'etale of degree $m(r,d,k,g)$,
and its universal kernel and quotient are vector bundles.
\end{proposition}
\begin{proof}
Put $B=\mathbb C[[t]]$, $K=\mathbb C((t))$, and
$N=\ell(\Quot_{d-e_*,s}(E))$.
A closed subset of a proper $B$-scheme disjoint from the closed
fibre is empty, since its image in $\operatorname{Spec}B$ is
closed and omits the closed point. We repeatedly use this fact.

\smallskip\noindent
\emph{Step 1: smoothing the pair.}
Local planarity gives a proper flat smoothing $\mathcal C$ over
$B$ \cite[Definition~5.9 and Proposition~5.10]{MRV15}.
The vanishings $H^2(C,\End E)=H^2(C,\mathcal O_C)=0$ give compatible
lifts of $E$ and of an ample line bundle to the infinitesimal
thickenings of $C$
\cite[Lemma~91.6.1, Tag~08VR]{Stk}.
Grothendieck existence algebraizes the lifts
\cite[Proposition~30.25.4, Tag~088C]{Stk}.
They are locally free along $C$ by completion and everywhere by
properness. The lifted line bundle is relatively ample
\cite[Lemma~30.21.4, Tag~0D2N]{Stk},
so the smoothing is projective.

Now fix any projective smoothing and bundle extension $\mathcal E$,
and choose a relatively very ample line bundle $\mathcal H$.
The smooth generic fibre is geometrically integral by semicontinuity
of $h^0(\mathcal O)$
\cite[Lemma~36.32.1, Tag~0BDN]{Stk}.
Constancy of Euler characteristics fixes the genus $g$, the bundle
degree $d$, and the fibre degree $h$ of $\mathcal H$
\cite[Lemma~36.32.2, Tag~0B9T]{Stk}.
The critical Hilbert polynomial is therefore
\[
 P_*(z)=shz+d-e_*+s(1-g).
\]

\smallskip\noindent
\emph{Step 2: finite \'etale transport.}
The full relative Quot
\[
 f:\mathcal Q=\Quot^{P_*}_{\mathcal C/B}(\mathcal E)
 \longrightarrow \operatorname{Spec}B
\]
is proper \cite[Lemma~108.6.3, Tag~0DPC]{Stk}
and has closed fibre $\Quot_{d-e_*,s}(E)$. At a point of this fibre,
$\Ext_C^1(F,Q)=0$ gives smoothness by quotient deformation theory
and the infinitesimal criterion
\cite[Remark~99.7.9, Tag~0CZU; Lemma~76.20.1, Tag~0APN]{Stk}.
The relative tangent space is $\Hom_C(F,Q)=0$
\cite[Remark~99.7.10, Tag~0CZV]{Stk},
so $f$ is \'etale there. Its non-\'etale locus is therefore empty
by properness. Thus $f$ is finite \'etale of degree $N$
\cite[Lemma~37.44.1, Tag~02LS]{Stk},
the closed fibre is reduced, and
\begin{equation}\label{eq:ade-critical-specialization}
 \#\mathcal Q_{\overline K}=N.
\end{equation}
The universal quotient is flat over $\mathcal Q$ and locally free
above the closed point by the fibrewise criterion
\cite[Lemma~37.16.8, Tag~0CZR]{Stk}.
Its nonfree locus on $\mathcal C\times_B\mathcal Q$ is empty by
properness over $B$, and its kernel is locally free by local splitting.

\smallskip\noindent
\emph{Step 3: passage to a smooth complex fibre.}
Holla's theorem is stated over $\mathbb C$, whereas our smoothing
has smooth generic fibre over $K=\mathbb C((t))$. To obtain a smooth
complex curve carrying the same count, we "spread out" the polarized
pair $(\mathcal C_K,\mathcal E_K,\mathcal H_K)$ and then take a
closed fibre of the resulting family. The projective embedding
and the bundle descend to a finitely generated
$\mathbb C$-subalgebra $A\subset K$, and $A$ can be enlarged so
that the curve family is smooth
\cite[Lemmas~32.10.1, 32.8.5, 32.8.9 and~32.10.3, Tags~01ZM, 01ZP, 0C0C and~0B8W]{Stk}.

Put $T=\operatorname{Spec}A$ and let $\eta$ be its generic point.
The map $\operatorname{Spec}K\to T$ factors through $\eta$, while
every closed point of $T$ has residue field $\mathbb C$. After
shrinking $T$, the curve fibres are geometrically integral by
semicontinuity of $h^0(\mathcal O)$. The genus and degrees are
constant, so the full relative Quot $\mathcal Q_T\to T$ with
Hilbert polynomial $P_*$ pulls back to $\mathcal Q_K$ over $K$.
By faithfully flat descent from $K$ to $\kappa(\eta)$, its generic
fibre is finite \'etale of degree $N$ and its universal quotient
is locally free there.

As in Step~2, the non-\'etale locus in $\mathcal Q_T$ and the
nonfree locus of the universal quotient in
$\mathcal C_T\times_T\mathcal Q_T$ are closed and disjoint from
the respective fibres over $\eta$. Both spaces are proper over
$T$, so the images of these loci are closed and omit $\eta$.
On their complement, $\mathcal Q_T$ is finite \'etale of degree
$N$, and the universal quotient and its kernel are vector bundles.
Taking a closed point of this open gives a smooth complex pair
$(D,V)$ whose full critical Quot consists of $N$ reduced points
with vector-bundle endpoints.

\smallskip\noindent
\emph{Step 4: comparison with a general bundle on $D$.}
The bundle $V$ obtained by specialization may lie outside Holla's required
general stable locus. We therefore vary the bundle on the fixed
curve $D$. Let $[V\to Q]$ be a point of $\Quot_{d-e_*,s}(V)$, with kernel $F$.
Since this Quot is finite and reduced, its tangent space
$\Hom_D(F,Q)$ vanishes, and then $\Ext_D^1(F,Q)=0$ because
$\chi(F^\vee\otimes Q)=D_{e_*}=0$. Applying Step~2's deformation argument on
smooth charts and removing the closed image of the non-\'etale
locus gives an open neighbourhood of $V$ in $\Bun_r^d(D)$
over which the full critical Quot is finite \'etale of degree $N$.
By Lemma~\ref{lem:ade-bun}, this open is connected and meets the
general stable locus of Holla's theorem
\cite[Proposition~4.1 and Theorem~4.2]{Hol04}. The degree there
is $m(r,d,k,g)$, so $N=m(r,d,k,g)$.
\end{proof}
\vspace{2mm}
\begin{remark}\label{rem:after-transport}
Proposition~\ref{prop:count-transport} reduces the count to the
geometry of the critical Quot scheme of a general bundle on $C$.
The next section gives two conditions,
\eqref{eq:quotient-boundary-hypothesis} and
\eqref{eq:critical-pair-nonspeciality}, under which this geometry
holds.
\end{remark}

\section{Geometry of the critical Quot scheme}
\label{sec:critical-geometry}

\begin{proposition}[A boundary criterion for the critical geometry]
\label{prop:boundary-generic-geometry}
Let $C$ be an integral projective complex curve of arithmetic genus
$g\ge2$. Fix $r,d$ and $0<k<r$, and put $s=r-k$ and
\[
 D_e=kd-re-ks(g-1).
\]
Write $\Coh^{\tf,\mathrm{nlf}}_{s,n}(C)$ for the closed locus of
non-locally-free sheaves in $\Coh^{\tf}_{s,n}(C)$. Assume that
\begin{equation*}\label{eq:quotient-boundary-hypothesis}
 \dim\Coh^{\tf,\mathrm{nlf}}_{s,d-e}(C)
 \le s^2(g-1)-1
 \qquad\text{whenever }D_e\le0.
 \tag{B}
\end{equation*}
Then the following hold.
\begin{enumerate}
\item\label{part:boundary-degree-bound}
On a nonempty open of $\Bun_r^d(C)$, every rank-$k$ subsheaf
$F\subset E$ satisfies $D_{\deg F}\ge0$.

\item\label{part:boundary-critical-finiteness}
Suppose that $e_*=(kd-ks(g-1))/r$ is an integer. After shrinking
this open, the full critical relative Quot morphism is finite,
possibly empty. Every geometric point of
$\Quot_{d-e_*,s}(E)$ is a sequence of vector bundles.

\item\label{part:boundary-pair-avoidance}
Under the integrality assumption in
\textup{(\ref{part:boundary-critical-finiteness})}, for every fixed
proper closed substack
\[
 Z\subsetneq\Bun_k^{e_*}(C)\times\Bun_s^{d-e_*}(C),
\]
a further nonempty open of ambient bundles has no critical endpoint
pair in $Z$. Moreover, for every fixed vector bundle $T$ of rank $k$
and degree $e_*$, a further nonempty open contains no critical
subsheaf isomorphic to $T$.

\item\label{part:boundary-nonspeciality}
Under the same integrality assumption, suppose, in addition, that
\begin{equation*}\label{eq:critical-pair-nonspeciality}
 \Hom_C(F,Q)=0
 \quad\text{for a general pair }
 (F,Q)\in\Bun_k^{e_*}(C)\times\Bun_s^{d-e_*}(C).
 \tag{N}
\end{equation*}
Then, on a nonempty open of $\Bun_r^d(C)$, the degree $e_*$ is
the maximal degree of a rank-$k$ subsheaf, and the full critical
Quot scheme is nonempty, finite and reduced. Every point is a
sequence of vector bundles satisfying
\[
 \Hom_C(F,Q)=\Ext_C^1(F,Q)=0.
\] In particular, the rank-$k$ subsheaves of $E$ of degree $e_*$ are
exactly the maximal rank-$k$ subbundles of $E$ in the sense of
Definition~\ref{def:maximal-subbundle}, and they are the points
of this Quot scheme.
\end{enumerate}
\end{proposition}

\begin{proof}
\noindent\emph{Common dimension estimates.}
Choose an irreducible finite-type smooth chart
$S\to\Bun_r^d(C)$ of constant relative dimension $m$, with
corresponding family $\mathcal E_S$ on $C\times S$, so that
$\dim S=r^2(g-1)+m$ by Lemma~\ref{lem:ade-bun}. Set
$H_e:=S\times_{\Bun_r^d(C)}\mathfrak X_e$, so that the square
\[
\begin{tikzcd}
 H_e \arrow[r] \arrow[d] & \mathfrak X_e \arrow[d, "\pi_E"] \\
 S \arrow[r] & \Bun_r^d(C)
\end{tikzcd}
\]
is cartesian. As noted after
Definition~\ref{def:saturated-quot-stack}, $H_e$ is the open
subscheme of the relative Quot scheme
$\Quot_{d-e,s}(\mathcal E_S/S)$ where the quotient is torsion-free.
It is of finite type by Lemma~\ref{lem:ade-bounded-charts}, and the
top arrow, being a base change of the chart map, is smooth of relative
dimension $m$.

Assume $D_e\le0$. The bundle locus of $\Coh^{\tf}_{s,d-e}(C)$ has
dimension $s^2(g-1)$ by Lemma~\ref{lem:ade-bun}, and its complement
has dimension at most $s^2(g-1)-1$
by~\eqref{eq:quotient-boundary-hypothesis}.
Corollary~\ref{cor:penalty-codim}, pulled back to $S$, therefore
gives
\[
 \dim H_e\le\dim S+D_e,
 \qquad
 \dim H_e^{\mathrm{nlf}}\le\dim S+D_e-1,
\]
where $H_e^{\mathrm{nlf}}$ is the locus with non-locally-free
quotient. We will construct nonempty good opens in $S$; their images
in $\Bun_r^d(C)$ are open, and the asserted properties hold there
because they are intrinsic to the ambient bundle.

We use repeatedly that if $X\to S$ is of finite type with
$\dim X<\dim S$, then its image is constructible by Chevalley's
theorem \cite[Theorem~29.23.3, Tag~054K]{Stk} and not dense, so the
complement of its closure is a nonempty open of $S$.

\smallskip\noindent
\emph{Proof of (\ref{part:boundary-degree-bound}): the degree bound.}
We have $D_e<0$ exactly when $e>(kd-ks(g-1))/r$, and
Lemma~\ref{lem:ade-bounded-charts} bounds above the degrees of
rank-$k$ subsheaves in the family $\mathcal E_S$. Hence only finitely
many degrees $e$ with $D_e<0$ occur on $S$, and for each of them
$\dim H_e\le\dim S+D_e<\dim S$. Since $S$ is irreducible, removing
the closures of the images of these finitely many $H_e$ leaves a
nonempty open $S^\circ\subset S$ over which every saturated rank-$k$
subsheaf satisfies $D_e\ge0$.

For an arbitrary rank-$k$ subsheaf $F\subset E$, the saturation has
degree\\ $\deg F+\ell$ with
$\ell=\operatorname{length}(F^{\mathrm{sat}}/F)\ge0$, so
\begin{equation}\label{eq:saturation-D}
 D_{\deg F}=D_{\deg F^{\mathrm{sat}}}+r\ell\ge0
\end{equation}
over $S^\circ$. The image of $S^\circ$ in $\Bun_r^d(C)$ is the
required open.

\smallskip\noindent
\emph{Proof of (\ref{part:boundary-critical-finiteness}):
finiteness and bundle endpoints.}
Replace $S$ by $S^\circ$, and assume that $e_*$ is an integer, so
that $D_{e_*}=0$ by \eqref{eq:critical-degree}. A rank-$k$ subsheaf
$F$ of degree $e_*$ is saturated: in \eqref{eq:saturation-D} the left
side vanishes and $D_{\deg F^{\mathrm{sat}}}\ge0$, so $\ell=0$. Thus
no critical quotient has torsion, at smooth or at singular points,
and the full critical relative Quot scheme equals $H_{e_*}$. By the common estimates,
$\dim H_{e_*}^{\mathrm{nlf}}\le\dim S-1$. Removing the closure of its
image leaves a nonempty open over which every critical quotient is a
vector bundle; its kernel is then a vector bundle by local splitting.

It remains to establish finiteness. Write $H=H_{e_*}$ over
the resulting open of $S$. The morphism $H\to S$ is proper
\cite[Lemma~108.6.3, Tag~0DPC]{Stk}, and
$\dim H\le\dim S$. Let $\eta$ be the generic point of $S$.
For every irreducible component $T$ of $H$ dominating $S$,
the dimension formula gives
\[
 \dim T_\eta=\dim T-\dim S\le0;
\]
since $T_\eta$ is nonempty, its dimension is zero.
Components not dominating $S$ have empty generic fibre.
Consequently $H_\eta$ is either empty or zero-dimensional.

If $H_\eta$ is empty, properness allows us to remove the
proper closed image of $H$. Otherwise, after shrinking $S$,
every fibre is zero-dimensional
\cite[Lemma~37.30.1, Tag~05F7]{Stk}.
In either case the remaining fibres are finite, since they
are of finite type over fields. Properness then makes
$H\to S$ finite
\cite[Lemma~37.44.1, Tag~02LS]{Stk}.
Taking the open image of this chart in $\Bun_r^d(C)$ proves
(\ref{part:boundary-critical-finiteness}); the critical Quot
may still be empty.

\smallskip\noindent
\emph{Proof of (\ref{part:boundary-pair-avoidance}):
endpoint incidences.}
Write $f:H\to S$ for the finite critical relative Quot.
Its universal endpoints $\mathcal F,\mathcal Q$ are vector bundles
by flatness and fibrewise local freeness
\cite[Lemma~37.16.8, Tag~0CZR]{Stk}. They define an endpoint morphism
\[
 a:H\longrightarrow
 \Bun_k^{e_*}(C)\times\Bun_s^{d-e_*}(C).
\]
For the prescribed proper closed substack $Z$, put $H_Z=a^{-1}(Z)$.
The ambient bundles admitting a critical pair in $Z$ form exactly
$f(H_Z)$. This image is closed because $H_Z\to S$ is finite.
We must show that it is a proper subset.

Since $H$ is of finite type, $a$ factors through a nonempty
finite-type open substack $\mathcal V$ of the product. By
Lemma~\ref{lem:ade-bun}, $\mathcal V$ is irreducible of dimension
$(k^2+s^2)(g-1)$. Since $Z$ is closed and does not contain the
dense open $\mathcal V$, the intersection $Z\cap\mathcal V$ is a
proper closed substack of $\mathcal V$, so
\[
 \dim(Z\cap\mathcal V)\le(k^2+s^2)(g-1)-1.
\]
Over a fixed endpoint pair, the extension groupoid is
$[\Ext_C^1(Q,F)/\Hom_C(Q,F)]$, of dimension
\[
 \extdim^1(Q,F)-\homdim(Q,F)
 =-\chi(Q^\vee\otimes F)
 =2ks(g-1).
\]
This formula holds in families: on endpoint charts, the relative
extension stack has a two-term presentation coming from the perfect
complex $Rp_*(Q^\vee\otimes F)$, whose formation commutes with
base change \cite[Lemma~36.30.4, Tag~0B91]{Stk}.
Consequently, after pullback to the ambient chart,
\[
 \dim H_Z
 \le (k^2+s^2)(g-1)-1+2ks(g-1)+m
 =r^2(g-1)+m-1
 =\dim S-1,
\]
since $r=k+s$ and $\dim S=r^2(g-1)+m$.
Thus $f(H_Z)$ is proper, and its complement is a nonempty open
over which every critical pair avoids $Z$.

For the fixed-bundle assertion, consider the framed incidence
\[
 I_T=
 \operatorname{Isom}_{C\times H/H}
       (T\boxtimes\mathcal O_H,\mathcal F).
\]
It is of finite type by
\cite[Lemma~108.4.1, Tag~0DLY]{Stk}, and its image in $S$
consists precisely of bundles admitting a critical kernel
isomorphic to $T$. Retaining the identification $F\simeq T$
leaves only $Q$ varying in the endpoint base. The same
extension-stack calculation therefore gives
\[
 \dim I_T
 \le s^2(g-1)+2ks(g-1)+m
 =\dim S-k^2(g-1)<\dim S.
\]
Removing the closure of the image of $I_T\to S$ therefore excludes
all critical kernels isomorphic to $T$, proving
(\ref{part:boundary-pair-avoidance}).

\smallskip\noindent
\emph{Proof of (\ref{part:boundary-nonspeciality}): nonemptiness, reducedness and maximality.}
Since $D_{e_*}=0$, Riemann--Roch gives
$\homdim(F,Q)=\extdim^1(F,Q)$ for all vector bundles $F,Q$ of ranks
$k,s$ and degrees $e_*,d-e_*$. Under
\eqref{eq:critical-pair-nonspeciality}, semicontinuity makes the
locus $\Hom_C(F,Q)\ne0$ a proper closed substack of the endpoint
product, and (\ref{part:boundary-pair-avoidance}) excludes it at
every critical point of a general ambient bundle. At these points
$\Hom_C(F,Q)=\Ext_C^1(F,Q)=0$; the first is the tangent space
\cite[Proposition~2.2.7]{HL}, so the critical Quot scheme is reduced.
It is nonempty for a general ambient bundle: for one pair $(F,Q)$
satisfying \eqref{eq:critical-pair-nonspeciality} we have
$\Ext_C^1(F,Q)=0$, so by Lemma~\ref{lem:split-smooth} the image of
$\pi_E$ contains a nonempty open substack of $\Bun_r^d(C)$, namely a
neighbourhood of $F\oplus Q$. By irreducibility it meets the open
constructed above.

Finally, every rank-$k$ subsheaf of $E$ has degree at most $e_*$ by
(\ref{part:boundary-degree-bound}), and those of degree $e_*$ are the
points of the critical Quot scheme, which are subbundles by
(\ref{part:boundary-critical-finiteness}). This proves
(\ref{part:boundary-nonspeciality}).
\end{proof}

\begin{corollary}[Stable bundles]\label{cor:stable-density}
Let $C$ be an integral projective complex curve of arithmetic genus $g\ge2$, let $a\ge1$ and $b$ be integers, and suppose that
\eqref{eq:quotient-boundary-hypothesis} holds for
$(r,d,k)=(a,b,j)$ for every $0<j<a$. Then stable bundles form a
nonempty dense open substack of $\Bun_a^b(C)$.
\end{corollary}
\begin{proof}
Every line bundle is stable, so let $a>1$. Stability is an open
condition \cite[Proposition~2.3.1]{HL}, and $\Bun_a^b(C)$ is
irreducible by Lemma~\ref{lem:ade-bun}, so it suffices to find one
nonempty open consisting of stable bundles. For $0<j<a$,
Proposition~\ref{prop:boundary-generic-geometry}
(\ref{part:boundary-degree-bound}) for $(r,d,k)=(a,b,j)$ gives a
nonempty open of $\Bun_a^b(C)$ on which every rank-$j$ subsheaf
$F\subset E$ satisfies
\[
 jb-a\deg F\ \ge\ j(a-j)(g-1)>0,
 \qquad\text{that is,}\qquad \mu(F)<\mu(E).
\]
The intersection of these opens over $0<j<a$ is nonempty and consists
of stable bundles.
\end{proof}
\begin{remark}[The two applications]\label{rem:after-critical-geometry}
The remaining sections verify
\eqref{eq:quotient-boundary-hypothesis} and
\eqref{eq:critical-pair-nonspeciality}: for integral locally planar
curves in rank two in Section~\ref{subsec:rank-two-planar}, and for
curves with one ADE singularity in all ranks in
Section~\ref{sec:ade-critical}, using the local estimates of
Section~\ref{sec:local-types}.
\end{remark}

\section{Locally planar curves in rank two}
\label{subsec:rank-two-planar}

Throughout this section, $C$ is an integral projective complex
curve of arithmetic genus $g\ge2$ with locally planar
singularities. We verify the two inputs of
Proposition~\ref{prop:boundary-generic-geometry} for $r=2$ and
$k=s=1$. Both follow from rank-one geometry.

\begin{lemma}[The rank-one boundary]
\label{lem:planar-rank-one-boundary}
For every integer $n$,
\[
 \dim\Coh^{\tf,\mathrm{nlf}}_{1,n}(C)\le g-2.
\]
In particular, condition~\eqref{eq:quotient-boundary-hypothesis}
holds in rank one.
\end{lemma}
\begin{proof}
Let $\bar J^n(C)$ be the degree-$n$ compactified Jacobian, which
parametrizes all rank-one torsion-free sheaves of degree $n$ because
$C$ is integral. By \cite[Theorem~A(ii)]{MRV15}, it has pure
dimension $g$ and its line-bundle locus is dense, so its
non-invertible locus has dimension at most $g-1$. Every rank-one
torsion-free sheaf on $C$ is simple, with automorphism group
$\mathbb G_m$, so the sheaf stack is a $\mathbb G_m$-gerbe over the
compactified Jacobian, and passing to it lowers the dimension by one.
\end{proof}

\begin{lemma}[Line-bundle nonspeciality]
\label{lem:planar-line-nonspeciality}
Let $a,b$ be integers with $b-a=g-1$. For a general pair of
line bundles $F,Q$ of degrees $a,b$, respectively,
\[
 \Hom_C(F,Q)=\Ext_C^1(F,Q)=0.
\]
\end{lemma}
\begin{proof}
Every line bundle of degree $g-1$ with a nonzero section is
$\mathcal O_C(D)$ for an effective Cartier divisor $D$ of
degree $g-1$, since $C$ is integral. Such a divisor is a point
of $\operatorname{Hilb}^{g-1}(C)$, namely the quotient
$\mathcal O_C\twoheadrightarrow\mathcal O_D$ with kernel
$\mathcal O_C(-D)$. Its tangent and obstruction spaces are
\[
 \Hom_C(\mathcal O_C(-D),\mathcal O_D)=H^0(C,\mathcal O_D(D)),
 \qquad
 \Ext_C^1(\mathcal O_C(-D),\mathcal O_D)=H^1(C,\mathcal O_D(D))=0.
\]
The second identity holds because $\mathcal O_D(D)$ has finite support. The
first has dimension $\operatorname{length}(D)=g-1$. Hence
$\operatorname{Hilb}^{g-1}(C)$ is smooth of dimension $g-1$ at
every Cartier divisor \cite[Proposition~2.2.8]{HL}.

Consequently, the image of the Abel map
$D\mapsto\mathcal O_C(D)$ on Cartier divisors has dimension at
most $g-1$.
This image is the locus of line bundles with a nonzero section,
which is closed by semicontinuity. Since
$\dim\Pic^{g-1}(C)=g$, a general line bundle $L$ of degree
$g-1$ has $H^0(C,L)=0$, and Riemann--Roch gives $H^1(C,L)=0$.
For any fixed $F$, tensoring by $F^{-1}$ identifies
$\Pic^b(C)$ with $\Pic^{g-1}(C)$. Applying this to
$L=F^{-1}\otimes Q$ proves the assertion.
\end{proof}

\begin{theorem}\label{thm:rank-two-planar}
Let $d\equiv g-1\pmod2$ and $e_*=(d-g+1)/2$.
There is a nonempty open substack of $\Bun_2^d(C)$, consisting
of stable bundles $E$, for which $e_*$ is the maximal degree
of a rank-one subsheaf and the full
$\Quot_{d-e_*,1}(E)$ is nonempty, finite and reduced.
Every point is a sequence of line bundles with
\[
 \Hom_C(F,Q)=\Ext_C^1(F,Q)=0.
\]
Consequently, every bundle in this open has exactly $2^g$
maximal line subbundles.
\end{theorem}
\begin{proof}
Lemma~\ref{lem:planar-rank-one-boundary} gives
condition~\eqref{eq:quotient-boundary-hypothesis}.
Since $d-2e_*=g-1$,
Lemma~\ref{lem:planar-line-nonspeciality} gives
condition~\eqref{eq:critical-pair-nonspeciality}.
Proposition~\ref{prop:boundary-generic-geometry}(\ref{part:boundary-nonspeciality}) therefore
supplies the asserted maximal degree and the nonempty finite
reduced full critical Quot with line-bundle endpoints and
the required vanishings.

By Lemma~\ref{lem:planar-rank-one-boundary} and Corollary~\ref{cor:stable-density}, we may shrink this open so that every $E$ in it is stable.
Finally, $C$ is locally planar and every $E$ in this open
satisfies the hypotheses of
Proposition~\ref{prop:count-transport}. Its critical Quot scheme
therefore has exactly $m(2,d,1,g)=2^g$ points
\cite[Corollary~4.3]{Hol04}, and these are the maximal line
subbundles of $E$.
\end{proof}

\section{Local types and their dimensions}\label{sec:local-types}

Throughout this section $C$ has one singular point $x$, without
any ADE assumption. Put
\[
 R=\widehat{\mathcal O}_{C,x},\qquad S=\overline R,
 \qquad\delta=\ell(S/R),\qquad\mathfrak c=(R:S),
\]
and let $K$ be the total quotient ring of $R$. Thus $S$ is a finite
product of power series rings $\mathbb C[[t]]$, one for each branch.
All modules are finitely generated, and $\ell$ is length over $R$,
which equals dimension over $\mathbb C$.

Let $M$ be a torsion-free $R$-module of constant rank $a$ on every
branch, so that $M\subset M\otimes_RK\simeq K^a$. Let $\bar M=SM$
be the $S$-submodule generated by $M$; equivalently
$\bar M=(S\otimes_RM)/\tors$. It is torsion-free of rank $a$ over
each factor of $S$, hence free, and $\mathfrak c\bar M=\mathfrak cM$.
After choosing a basis of $\bar M$,
\begin{equation}\label{eq:normalized-lattice}
 \mathfrak cS^a\subset M\subset\bar M=S^a\subset K^a .
\end{equation}
Every endomorphism of $M$ extends uniquely to $K^a$ and preserves
$\bar M$, so $\End_R(M)$ is the set of matrices in
$\operatorname{Mat}_a(K)$ preserving $M$, and
\begin{equation}\label{eq:end-chain}
 \mathfrak c\operatorname{Mat}_a(S)\subset\End_R(M)
 \subset\operatorname{Mat}_a(S)=\End_S(\bar M).
\end{equation}
These inclusions have finite colength, and we define
\begin{equation}\label{eq:ade-local-defect}
 c_R(M)=a^2\delta-
 \ell\bigl(\End_S(\bar M)/\End_R(M)\bigr).
\end{equation}
A local type means the isomorphism class of the completed stalk.
Every rank-$a$ torsion-free sheaf has constant completed branch
rank $a$.

\subsection{The local Gorenstein invariant}

\begin{definition}\label{def:stable-end}
Write $(-)^*=\Hom_R(-,R)$. An endomorphism of $M$ is
\emph{stably trivial} if it factors through a finitely generated
free module; equivalently, if it lies in the image $\mathcal P$ of
\[
 M^*\otimes_RM\longrightarrow\End_R(M),\qquad
 f\otimes m\longmapsto\bigl(x\mapsto f(x)m\bigr).
\]
This image is a two-sided ideal, and the
\emph{stable endomorphism ring} of $M$ is
$\underline{\End}_R(M)=\End_R(M)/\mathcal P$.
\end{definition}

In this subsection $R$ is Gorenstein. We use two standard facts.
\begin{enumerate}
\item\label{it:reflexive} Every torsion-free module $N$ is
 reflexive, and $\Ext^1_R(N,R)=0$
 \cite[Theorems~3.3.7 and~3.3.10(d)]{BH98}; see also
 \cite[\S2.1 and Lemma~2.5.3]{Coo93}.
\item\label{it:length} If $T$ has finite length, then
 $\ell(\Ext^1_R(T,R))=\ell(T)$. Indeed $\Ext^i_R(T,R)=0$ for
 $i\ne1$ \cite[Theorem~3.3.10(c)]{BH98}, so $\Ext^1_R(-,R)$ is
 exact on finite-length modules, and
 $\Ext^1_R(\mathbb C,R)\simeq\mathbb C$ because the canonical
 module $R$ has type one \cite[Definition~3.3.1]{BH98}.
\end{enumerate}

\begin{lemma}[Trace duality]\label{lem:trace-duality}
For a submodule $\Lambda\subset\operatorname{Mat}_a(K)$
spanning $\operatorname{Mat}_a(K)$ over $K$, put
\[
 \Lambda^\vee=\{\psi\in\operatorname{Mat}_a(K):
 \operatorname{tr}(\psi\lambda)\in R
 \text{ for all }\lambda\in\Lambda\}.
\]
Then $\Lambda^{\vee\vee}=\Lambda$, and
$\ell(\Lambda'^\vee/\Lambda^\vee)=\ell(\Lambda/\Lambda')$ for two
such submodules $\Lambda'\subset\Lambda$. In particular
$\ell(S/\mathfrak c)=2\delta$.
\end{lemma}
\begin{proof}
The trace form is perfect over $K$, so it identifies
$\Lambda^\vee$ with $\Lambda^*$, and
\eqref{it:reflexive} gives $\Lambda^{\vee\vee}=\Lambda$.
Dualizing $0\to\Lambda'\to\Lambda\to\Lambda/\Lambda'\to0$ gives
\[
 0\to\Lambda^*\to\Lambda'^*\to
 \Ext^1_R(\Lambda/\Lambda',R)\to0
\]
by \eqref{it:reflexive}, and \eqref{it:length} gives the length.
For $a=1$ we have $R^\vee=R$ and $S^\vee=\mathfrak c$, so
$\ell(R/\mathfrak c)=\ell(S/R)=\delta$.
\end{proof}
\begin{proposition}\label{prop:local-gorenstein-penalty}
Let $R$ be Gorenstein and $M$ a torsion-free $R$-module of
constant branch rank. Then
\begin{equation}\label{eq:ade-stable-end}
 2c_R(M)=\ell\bigl(\underline{\End}_R(M)\bigr),
\end{equation}
where $\underline{\End}_R(M)$ is the stable endomorphism ring of
\cref{def:stable-end}. In particular $c_R(M)\ge0$, with equality
precisely when $M$ is free.
\end{proposition}
\begin{proof}
Choose a basis of $\bar M$, so that $\bar M=S^a$ and all modules
below are lattices in $\operatorname{Mat}_a(K)$. The main observation
is that trace duality carries the chain
\begin{equation}\label{eq:self-dual-chain}
 \mathfrak c\operatorname{Mat}_a(S)\subset\mathcal P
 \subset\End_R(M)\subset\operatorname{Mat}_a(S)=\End_S(\bar M)
\end{equation}
to itself, interchanging the two inner terms and the two outer terms.

By \cref{def:stable-end}, $\mathcal P$ is spanned by the matrices
$mf$, where $m\in M$ is a column vector and $f\in M^*$ is a row
vector; it spans $\operatorname{Mat}_a(K)$, since $M$ spans $K^a$
and $M^*$ spans its dual. For $\psi\in\operatorname{Mat}_a(K)$ we
have $\operatorname{tr}(\psi\,mf)=f(\psi m)$, so
$\psi\in\mathcal P^\vee$ if and only if $\psi M\subset M^{**}$.
Since $M^{**}=M$ by \eqref{it:reflexive},
\[
 \mathcal P^\vee=\End_R(M),\qquad
 \End_R(M)^\vee=\mathcal P,
\]
the second by Lemma~\ref{lem:trace-duality}. Using
$S^\vee=\mathfrak c$, we have
$\operatorname{Mat}_a(S)^\vee=\mathfrak c\operatorname{Mat}_a(S)$.
Dualizing $\End_R(M)\subset\operatorname{Mat}_a(S)$ therefore gives
the first inclusion in \eqref{eq:self-dual-chain}, and by
Lemma~\ref{lem:trace-duality} and \eqref{eq:ade-local-defect},
\[
 \ell\bigl(\mathcal P/\mathfrak c\operatorname{Mat}_a(S)\bigr)
 =\ell\bigl(\operatorname{Mat}_a(S)/\End_R(M)\bigr)
 =a^2\delta-c_R(M).
\]
Since \eqref{eq:self-dual-chain} has total length
$a^2\ell(S/\mathfrak c)=2a^2\delta$, again by
Lemma~\ref{lem:trace-duality}, we obtain
\[
 \ell\bigl(\underline{\End}_R(M)\bigr)
 =\ell\bigl(\End_R(M)/\mathcal P\bigr)
 =2a^2\delta-2\bigl(a^2\delta-c_R(M)\bigr)=2c_R(M).
\]

Consequently $c_R(M)\ge0$, with equality if and only if
$\mathrm{id}_M\in\mathcal P$. In that case $\mathrm{id}_M$ factors
through a finitely generated free module, so $M$ is a direct summand
of it and hence free, because $R$ is local. Conversely, every
endomorphism of a free module is stably trivial.
\end{proof}

The statements of this subsection, and their proofs, are unchanged
when $\mathbb C$ is replaced by an algebraically closed extension
field.
\subsection{The fixed-type estimate}

We bound the dimension of the locus of torsion-free sheaves with a
prescribed local type, using the conductor description of such
sheaves \cite[\S\S2.2--2.3]{Coo93}.

\begin{lemma}[Fixed-type closure estimate]
\label{prop:ade-type-dimension}
Fix $a,b$ and a torsion-free $R$-module $M$ of constant branch
rank $a$. The locus $\mathfrak T_M(a,b)$ of sheaves with
completed stalk isomorphic to $M$ is locally constructible in
$\Coh^{\tf}_{a,b}(C)$. Its reduced closure satisfies
\begin{equation}\label{eq:ade-type-dimension}
 \dim\overline{\mathfrak T_M(a,b)}
 \le a^2(g-1)-c_R(M).
\end{equation}
The same assertion holds on finite-type smooth charts and after
algebraically closed field extension.
\end{lemma}

\begin{proof}
\noindent\emph{Step 1: gluing data.}
Let $\iota:Z\hookrightarrow C$ be the conductor subscheme and
$\widetilde Z=\nu^{-1}(Z)$. Let $B$ be a torsion-free sheaf of rank
$a$ on $C$ and put $V=B^\nu$, so that $B\subset\nu_*V$. Since
$\mathfrak c\cdot\nu_*V=\mathfrak c\cdot B$ by
\eqref{eq:normalized-lattice}, the restriction $W=B|_Z$ is an
$\mathcal O_Z$-submodule of $(\nu_*V)|_Z$, and $B$ is its preimage:
\begin{equation}\label{eq:glued-sheaf}
 B=\ker\bigl(\nu_*V\longrightarrow
 \iota_*\bigl((\nu_*V)|_Z/W\bigr)\bigr).
\end{equation}
Thus $B$ is determined by the pair $(V,W)$.

\smallskip\noindent\emph{Step 2: gluing data of local type $M$.}
Since $Z$ is supported at $x$, we have $\mathcal O_Z=R/\mathfrak c$
and $\nu_*\mathcal O_{\widetilde Z}=S/\mathfrak c$. An $S$-basis of
the completion of $\nu_*V$ at $x$ identifies $(\nu_*V)|_Z$ with
$(S/\mathfrak c)^a$, $W$ with an $R/\mathfrak c$-submodule of it,
and the completed stalk $\widehat B$ with the preimage of $W$ in
$S^a$. The group $G=\operatorname{GL}_a(S/\mathfrak c)$ acts on the
set of $R/\mathfrak c$-submodules of $(S/\mathfrak c)^a$, and a
change of basis moves $W$ by an element of $G$. Fix an embedding
$M\subset S^a$ with $SM=S^a$, so that $\mathfrak cS^a\subset M$, and
let $O_M$ be the $G$-orbit of $M/\mathfrak cS^a$. Then
\begin{equation}\label{eq:orbit-criterion}
 \widehat B\simeq M\iff W\in O_M .
\end{equation}
Indeed, $S\widehat B=S^a$, so an isomorphism $M\to\widehat B$ is
given by some $\tilde g\in\operatorname{GL}_a(K)$ with
$\tilde gS^a=S\,\tilde gM=S^a$; thus
$\tilde g\in\operatorname{GL}_a(S)$, and its reduction carries
$M/\mathfrak cS^a$ to $W$; compare
\cite[Proposition~2.2.1 and Lemma~2.2.2]{Coo93}. Conversely, every
element of $G$ lifts to $\operatorname{GL}_a(S)$, because
$\mathfrak c$ lies in the radical of $S$, and the lift carries $M$
onto the preimage of $W$. In particular the condition $W\in O_M$
depends neither on the basis nor on the embedding of $M$.

\smallskip\noindent\emph{Step 3: dimension of the orbit.}
The stabilizer of $M/\mathfrak cS^a$ in $G$ is the unit group of the
finite-dimensional algebra
$\End_R(M)/\mathfrak c\operatorname{Mat}_a(S)$. Hence $O_M$ is a
smooth locally closed subvariety of the Grassmannian of subspaces of
codimension $\lambda=\ell(S^a/M)$ in $(S/\mathfrak c)^a$, and by
\eqref{eq:ade-local-defect}, without any Gorenstein hypothesis,
\[
 \dim O_M
 =a^2\ell(S/\mathfrak c)
  -\ell\bigl(\End_R(M)/\mathfrak c\operatorname{Mat}_a(S)\bigr)
 =a^2\delta-c_R(M).
\]

\smallskip\noindent\emph{Step 4: the parameter stack.}
By Steps~1 and~2, a torsion-free sheaf $B$ of rank $a$, degree $b$
and local type $M$ is the same as a pair $(V,W)$, where $V$ is a
vector bundle of rank $a$ on $Y$ and $W\subset(\nu_*V)|_Z$ lies in
$O_M$. Here $\ell(\nu_*V/B)=\lambda$, so $V$ has degree
$b-a\delta+\lambda$ by \eqref{eq:ade-normalization-degree}. Let
$\mathfrak P_M$ be the stack of such pairs, and let
\[
 \mathfrak P_M\longrightarrow\Bun_a^{\,b-a\delta+\lambda}(Y),
 \qquad(V,W)\longmapsto V,
\]
forgetting the gluing datum. The $S/\mathfrak c$-bases of $(\nu_*V)|_Z$
form a $G$-torsor over the base, and $\mathfrak P_M$ is the
associated bundle with fibre $O_M$: each basis identifies the fibre
over $V$ with $O_M$. Hence this morphism is smooth of relative
dimension $\dim O_M$, and since $g(Y)=g-\delta$,
Lemma~\ref{lem:ade-bun} and Step~3 give
\[
 \dim\mathfrak P_M=a^2(g-\delta-1)+a^2\delta-c_R(M)
 =a^2(g-1)-c_R(M).
\]

\smallskip\noindent\emph{Step 5: the morphism to torsion-free sheaves.}
In a family of pairs $(V,W)$, the quotient $(\nu_*V)|_Z/W$ in
\eqref{eq:glued-sheaf} is locally free of rank $\lambda$ over the
base, being pulled back from the universal quotient on the
Grassmannian containing $O_M$. Hence the kernel $B$ is flat and
compatible with base change, and $(V,W)\mapsto B$ defines a morphism
\[
 \Phi:\mathfrak P_M\longrightarrow\Coh^{\tf}_{a,b}(C),
\]
whose image is exactly $\mathfrak T_M(a,b)$ by Step~4. Since $\nu^*B\to V$ is
surjective, an automorphism of $(V,W)$ inducing the identity on $B$
is the identity. Thus $\Phi$ is faithful, hence representable by
algebraic spaces \cite[Lemma~94.15.2, Tag~04Y5]{Stk}.

\smallskip\noindent\emph{Step 6: dimension of the image.}
Let $U\to\Coh^{\tf}_{a,b}(C)$ be a finite-type affine smooth chart of
relative dimension $m$, and let $B_u$ denote the fibres of its
universal sheaf. This sheaf is a quotient of a finite sum of powers
of a very ample line bundle \cite[Lemma~30.16.1(2), Tag~0B5T]{Stk},
so after a common twist all $B_u$ are globally generated. Then so
are the $B_u^\nu$, as quotients of $\nu^*B_u$, and a further common
twist kills their $H^1$ by Serre duality. Hence the $B_u^\nu$ lie in
a finite-type open substack of the bundle stack of $Y$, presented by
a Quot scheme with fixed Hilbert polynomial
\cite[Lemma~108.6.4, Tag~0DPD]{Stk}; let
$\mathfrak P^0\subset\mathfrak P_M$ be its preimage. By Step~5 and
\cite[Lemma~108.4.1, Tag~0DLY]{Stk},
$P_U=\mathfrak P^0\times_{\Coh^{\tf}_{a,b}(C)}U$ is a finite-type
algebraic space, smooth of relative dimension $m$ over
$\mathfrak P^0$, with image $\mathfrak T_M(a,b)\cap U$ in $U$. This
image is constructible by Chevalley's theorem
\cite[Theorem~67.29.2, Tag~0ECX]{Stk}, and by Step~4
\[
 \dim\overline{\mathfrak T_M(a,b)\cap U}
 \le\dim P_U\le a^2(g-1)-c_R(M)+m.
\]
Subtracting $m$ gives \eqref{eq:ade-type-dimension}.
\end{proof}

\section{One ADE singularity}\label{sec:ade-critical}

\begin{definition}\label{def:ade-singularity}
A singular point $x$ of a complex curve $C$ is an
\emph{ADE singularity} if
\[
 \widehat{\mathcal O}_{C,x}
 \simeq \mathbb C[[u,v]]/(f),
\]
where, after a formal change of coordinates, $f$ has one of the
following forms:
\[
\begin{array}{ccl}
 A_n &:& v^2-u^{n+1}, \qquad n\geq1,\\
 D_n &:& uv^2+u^{n-1}, \qquad n\geq4,\\
 E_6 &:& v^3+u^4,\\
 E_7 &:& v^3+u^3v,\\
 E_8 &:& v^3+u^5.
\end{array}
\]
These are the simple plane curve singularities
\cite[p.~417]{GK85} (our sign convention is different, but over complex numbers up to change of coordinates they give the same thing). In particular, their completed local rings
are reduced, Gorenstein, and of finite Cohen--Macaulay type:
there are only finitely many isomorphism classes of
indecomposable finitely generated torsion-free modules
\cite[Satz~1, p.~418, and Korollar, p.~419]{GK85}.
\end{definition}

Throughout this section $C$ has arithmetic genus $g\ge2$ and
exactly one singular point $x$, of ADE type. We use the notation
$R,S,\delta,c_R$ of Section~\ref{sec:local-types}; for a sheaf $B$,
write $\widehat B$ for its completed stalk at $x$.

We verify \eqref{eq:quotient-boundary-hypothesis} and
\eqref{eq:critical-pair-nonspeciality} in all ranks. In
Section~\ref{subsec:rank-two-planar} the quotients had rank one, so
both conditions were statements about the compactified Jacobian and
about line bundles. For quotients of rank $s\ge2$ we instead stratify
the non-locally-free sheaves by their type at $x$: there are finitely
many types (Lemma~\ref{lem:ade-types}), each of positive codimension
(Proposition~\ref{prop:local-gorenstein-penalty} and
Lemma~\ref{prop:ade-type-dimension}). Nonspeciality is then no longer
elementary and follows the argument of Hirschowitz and Bhosle
\cite{Hir88,Bho06}.

\subsection{The boundary condition}

\begin{lemma}[Finitely many types in fixed rank]\label{lem:ade-types}
Let $a\ge1$. Up to isomorphism there are only finitely many
torsion-free $R$-modules of rank $a$ on every branch. In particular,
only finitely many local types occur among the completed stalks at
$x$ of rank-$a$ torsion-free sheaves on $C$.
\end{lemma}

\begin{proof}
By \cite[Satz~1, p.~418, and Korollar, p.~419]{GK85}, $R$ has
finitely many indecomposable torsion-free modules $N_1,\dots,N_p$,
some of which may vanish on some branches. By Krull--Schmidt, every
torsion-free module is $\bigoplus_iN_i^{n_i}$. If it has rank $a$ on
each of the $q$ branches, then summing the ranks over the branches
gives $n_i\le aq$, because each $N_i$ has positive rank on at least
one branch. This leaves finitely many modules. The last assertion
holds because completed stalks of rank-$a$ torsion-free sheaves have
rank $a$ on every branch (Section~\ref{sec:local-types}).
\end{proof}

\begin{corollary}[The ADE boundary bound]\label{cor:ade-boundary}
For every $a\ge1$ and every integer $b$,
\[
 \dim\Coh^{\tf,\mathrm{nlf}}_{a,b}(C)
 \le a^2(g-1)-1.
\]
Thus condition~\eqref{eq:quotient-boundary-hypothesis} holds in
every rank and degree.
\end{corollary}
\begin{proof}
By Lemma~\ref{lem:ade-types}, the nonfree locus is the union of
finitely many nonfree type loci. For each such type $M$,
Lemma~\ref{prop:ade-type-dimension} and
Proposition~\ref{prop:local-gorenstein-penalty} give
\[
 \dim\overline{\mathfrak T_M(a,b)}
 \le a^2(g-1)-c_R(M)\le a^2(g-1)-1,
\]
since $c_R(M)$ is a positive integer. Taking the finite union,
with dimensions interpreted on finite-type smooth charts,
proves the assertion.
\end{proof}

\begin{remark}
Here the two local hypotheses have distinct roles. Gorenstein
duality makes the penalty positive for each nonfree type;
finite Cohen--Macaulay type allows these bounds to control the
whole nonfree locus. Positivity for individual types alone does
not bound a union of moving types.
\end{remark}

\begin{corollary}\label{cor:ade-stable-density}
Stable bundles form a nonempty dense open of $\Bun_a^b(C)$.
\end{corollary}
\begin{proof}
This follows from \cref{cor:ade-boundary} and \cref{cor:stable-density}.
\end{proof}

\subsection{Nonspeciality}

\begin{remark}
The proof below uses parts
(\ref{part:boundary-degree-bound})--(\ref{part:boundary-pair-avoidance})
of Proposition~\ref{prop:boundary-generic-geometry}. These require
only \eqref{eq:quotient-boundary-hypothesis}, which holds by
Corollary~\ref{cor:ade-boundary}.
\end{remark}
\vspace{2mm}
\begin{proposition}[Generic nonspeciality]
\label{prop:ade-nonspeciality}
For general stable bundles $U,V$ of any fixed positive ranks
and degrees,
\[
\begin{array}{ll}
H^0(C,U^\vee\otimes V)=0,
  &\mu(V)-\mu(U)\le g-1,\\
H^1(C,U^\vee\otimes V)=0,
  &\mu(V)-\mu(U)\ge g-1.
\end{array}
\]
\end{proposition}
\begin{proof}
The argument is that of Hirschowitz \cite{Hir88}, as carried out for
nodal curves in \cite[Proposition~2.9, Lemma~2.10 and proof of
Theorem~1]{Bho06}. We indicate its inputs and why they hold on $C$.

The vanishings are open conditions on the irreducible product of
bundle stacks, in which stable pairs are dense
(Lemma~\ref{lem:ade-bun}, Corollary~\ref{cor:ade-stable-density}), so
it suffices to exhibit one pair of bundles with the required
vanishing; this replaces \cite[Lemma~2.4]{Bho06}. Since $C$ is
Gorenstein, Serre duality gives
$H^1(C,U^\vee\otimes V)^\vee\simeq\Hom_C(V,U\otimes\omega_C)$, and
$U\mapsto U\otimes\omega_C$ is an isomorphism of bundle stacks; hence
the $H^1$ assertion for $(U,V)$ is the $H^0$ assertion for
$(V,U\otimes\omega_C)$, and we prove the latter. The inputs needed
are the following.
\begin{enumerate}
\item[(a)] \emph{Compositions vanish.} For a general pair $(U,V)$,
every $f:U\to V$ and $h:V\to U\otimes\omega_C$ satisfy $h\circ f=0$.
This is \cite[Proposition~2.9(2)]{Bho06}, whose proof uses only that
the bundle stacks are smooth and irreducible and that $\omega_C$ is
invertible.
\item[(b)] \emph{Subsheaves of a general bundle.} For a general
bundle, every subsheaf of intermediate rank satisfies the degree
bound; at the critical degree there are finitely many, with
vector-bundle endpoints; and no fixed bundle occurs among them. For
nodal curves this is \cite[Proposition~2.8(1)--(3)]{Bho06}, proved
there from the local structure at a node. Here it is
Proposition~\ref{prop:boundary-generic-geometry}
(\ref{part:boundary-degree-bound})--(\ref{part:boundary-pair-avoidance}),
which applies by Corollary~\ref{cor:ade-boundary}.
\end{enumerate}

Granting (a) and (b), the proofs of \cite[Lemma~2.10]{Bho06}, for
equal ranks and $\deg V=\deg U+\rk U\,(g-1)$, and of
\cite[Theorem~1]{Bho06}, for the remaining ranks and degrees, apply
without change, with two additional observations. First, in Lemma~2.10 the bundle $V$
must avoid the finitely many critical quotients of $U$, which depend
on $U$. This holds for a general pair: the relative critical Quot
scheme is finite over charts of the first factor, so the exceptional
pairs form a constructible set that is nowhere dense in each fibre
\cite[Lemma~108.4.1, Tag~0DLY]{Stk}. Second, where the proof of
Theorem~1 passes to a subsheaf of a general bundle $G$ of smaller
degree, we take the kernel of a surjection from $G$ onto a torsion
sheaf supported at smooth points. It is a vector bundle, so
\cite[Lemma~2.3]{Bho06} is not needed.
\end{proof}

\subsection{The critical Quot and its count}

\begin{theorem}\label{thm:ade-critical-quot}
Assume the critical congruence \eqref{eq:critical-congruence}. There is a nonempty open
substack of stable rank-$r$, degree-$d$ bundles $E$ such that
$e_*$ is the maximal degree of a rank-$k$ subsheaf and the full
scheme $\Quot_{d-e_*,r-k}(E)$ is nonempty, finite and reduced.
Every point is a sequence of stable vector bundles and satisfies
\[
 \Hom_C(F,Q)=\Ext_C^1(F,Q)=0.
\]
\end{theorem}
\begin{proof}
Corollary~\ref{cor:ade-boundary} gives
\eqref{eq:quotient-boundary-hypothesis}. At the critical degree,
$\mu(Q)-\mu(F)=g-1$, so Proposition~\ref{prop:ade-nonspeciality}
gives \eqref{eq:critical-pair-nonspeciality}.
Proposition~\ref{prop:boundary-generic-geometry}
(\ref{part:boundary-nonspeciality}) supplies the critical Quot
geometry, maximality and the two vanishings.
Intersect its open with the dense stable locus of
Corollary~\ref{cor:ade-stable-density}. Finally, apply
(\ref{part:boundary-pair-avoidance}) to the proper closed
endpoint locus where either bundle is not stable, again using
Corollary~\ref{cor:ade-stable-density}. This makes every critical
endpoint stable after a further shrinking.
\end{proof}

\begin{theorem}\label{thm:ade-counting}
Assume the critical congruence \eqref{eq:critical-congruence}. A general stable rank-$r$,
degree-$d$ bundle on $C$ has exactly $m(r,d,k,g)$ maximal rank-$k$
subbundles, where $m(r,d,k,g)$ is given by
\eqref{eq:ade-holla-count}.
\end{theorem}
\begin{proof}
ADE singularities are locally planar. For the bundles of
Theorem~\ref{thm:ade-critical-quot}, the maximal rank-$k$
subbundles are the points of the full critical Quot scheme by
Proposition~\ref{prop:boundary-generic-geometry}
(\ref{part:boundary-nonspeciality}), and
Proposition~\ref{prop:count-transport} shows that this scheme is
reduced of length $m(r,d,k,g)$.
\end{proof}

\subsection{The degree threshold and the refined bound}

\begin{remark}[Beyond the critical count]
\label{rem:ade-beyond-counting}
The critical counting problem for ADE curves has already been
settled by conditions (B) and (N), together with the transport
theorem. We now record two further consequences of the same
methods. First, we recover Hirschowitz's maximal-subsheaf
degree formula without imposing the critical congruence.
Second, we refine the boundary estimate by relating the
expected Quot dimension to the completed local module type
of a torsion-free quotient. In particular, the penalty
$c_R(\widehat Q)$ constrains which non-locally-free quotient
types can occur for a general ambient bundle.
\end{remark}

\begin{corollary}[The generic threshold]\label{cor:ade-hirschowitz}
For general stable $E$ of rank $r$ and degree $d$, the maximal
degree of a rank-$k$ subsheaf is
\[
 e_0=\left\lfloor\frac{kd-k(r-k)(g-1)}r\right\rfloor,
\]
and it is attained by a vector subbundle with stable endpoints.
Equivalently,
\begin{equation}\label{eq:Hirschowitz-counting}
 s_k(E)=kd-r e_{\max}(E)=k(r-k)(g-1)+\epsilon,
 \qquad 0\le\epsilon<r,
\end{equation}
where $s_k(E)\equiv kd\pmod r$ determines the remainder.
\end{corollary}
\begin{proof}
Put $\epsilon=D_{e_0}$, so that $0\le\epsilon<r$. For general stable
bundles $F,Q$ of ranks $k,s$ and degrees $e_0,d-e_0$,
\[
 \mu(Q)-\mu(F)=g-1+\frac{\epsilon}{ks}\ge g-1,
\]
so $\Ext_C^1(F,Q)=0$ by Proposition~\ref{prop:ade-nonspeciality}.
The sequences in $\mathfrak X_{e_0}$ whose endpoints are stable
vector bundles form an open substack
(Corollary~\ref{cor:ade-stable-density}), which contains the split
sequence $0\to F\to F\oplus Q\to Q\to0$. By
Lemma~\ref{lem:split-smooth} its image contains a nonempty open
substack of $\Bun_r^d(C)$, whose bundles therefore have a subbundle
of degree $e_0$ with stable endpoints, and which meets the dense
stable locus. No larger degree occurs, by
Proposition~\ref{prop:boundary-generic-geometry}
(\ref{part:boundary-degree-bound}), since $D_{e_0+1}=\epsilon-r<0$.
\end{proof}
\begin{remark}[Comparison with the critical case]
\label{rem:threshold-vs-critical}
When $e_*$ is an integer, $e_0=e_*$ and
Theorem~\ref{thm:ade-critical-quot} says much more: the full Quot
scheme at the maximal degree is nonempty, finite and reduced, and
every one of its points is a sequence of vector bundles.
Corollary~\ref{cor:ade-hirschowitz} drops the congruence and retains
only the value of the maximal degree and the existence of one
maximal subbundle; the Quot scheme at degree $e_0$ has expected
dimension $\epsilon=D_{e_0}$, so for $\epsilon>0$ it need not be
finite, and its points need not all be sequences of vector bundles.
\end{remark}

\begin{proposition}[Type-dependent refinement]
\label{prop:ade-generic-geometry}
Fix $r,d$ and $0<k<r$, and put $s=r-k$.
\begin{enumerate}
\item
For a completed quotient type $N$ of constant branch rank $s$,
the locus $\mathfrak X_{e,N}$ of saturated sequences with
$\widehat Q\simeq N$ is locally constructible and satisfies
\begin{equation}\label{eq:ade-extension-penalty}
 \dim\overline{\mathfrak X_{e,N}}
 \le r^2(g-1)+D_e-c_R(N).
\end{equation}
\item\label{part:ade-lower-bound}
On a nonempty open of $\Bun_r^d(C)$, every saturated sequence
of any degree $e$ satisfies
\begin{equation}\label{eq:ade-generic-bound}
 D_e\ge c_R(\widehat Q).
\end{equation}
\end{enumerate}
\end{proposition}
\begin{proof}
The first assertion follows from
Lemma~\ref{prop:ade-type-dimension} and
Corollary~\ref{cor:penalty-codim}, with constructibility and
closures interpreted on finite-type charts.
For the second, argue on a chart $S$ as in the proof of
Proposition~\ref{prop:boundary-generic-geometry}
(\ref{part:boundary-degree-bound}). By Lemma~\ref{lem:ade-types} the
penalties of rank-$s$ types have a finite maximum, so a violation
$D_e<c_R(N)$ bounds $e$ below, while
Lemma~\ref{lem:ade-bounded-charts} bounds it above. Hence only
finitely many pairs $(e,N)$ can violate the bound, and for each
\eqref{eq:ade-extension-penalty} makes the image in $S$ not dense.
Removing the closures of these images proves the assertion.
\end{proof}

\appendix
\section{The explicit \texorpdfstring{$A_n$}{A\_n} penalty}\label{app:an-example}

\begin{remark}
The counting theorems use only that $c_R(M)>0$ for nonfree $M$.
For $A_n$ singularities we compute $c_R$ in terms of the standard
module multiplicities (Corollary~\ref{cor:an-explicit-penalty}),
which turns~\eqref{eq:ade-generic-bound} into an explicit numerical
bound.
\end{remark}
\vspace{2mm}
Fix $n\ge1$ and put
\[
 R=\mathbb C[[u,v]]/(v^2-u^{n+1}),\qquad
 A=\mathbb C[[u]],\qquad
 \delta=\lfloor(n+1)/2\rfloor.
\]
Its normalization and conductor are
\[
 S=A\oplus A(v/u^\delta),\qquad
 \mathfrak c=u^\delta S,
 \qquad \ell(S/R)=\delta.
\]
Set $I_i=(u^i,v)$ for $0\le i\le\delta$, so that $I_0=R$.
For a torsion-free $R$-module $M$ of constant branch rank $a$,
the classification in
\cite[Appendix, (1.1), (2.0)--(2.1), pp.~423--424]{GK85}
gives a unique balanced expression
\[
 M\simeq\bigoplus_{i=0}^{\delta}I_i^{\oplus b_i},
 \qquad \sum_{i=0}^{\delta}b_i=a.
\]
For odd $n=2\delta-1$, the two branch-supported indecomposables
occur with equal multiplicities, because the total branch ranks
are equal. Pairing them gives the block $I_\delta\simeq S$;
this block is not indecomposable. For even $n$, all the displayed
$I_i$ are indecomposable.

We record the normalized Hom calculation; compare the classical
pairwise computations in
\cite[proof of Lemma~2.5.9, pp.~31--32]{Coo93}.
As an $A$-module $R=A\oplus Av$, and since $v^2=u^{n+1}$ and
$i\le\delta\le n+1$, we have $I_i=u^iA\oplus Av$. Put
\[
 J_i=u^{-i}I_i=A\oplus u^{-i}Av\subset S=A\oplus u^{-\delta}Av .
\]
Since $1\in J_i$ and $u^iJ_i\subset R$, a homomorphism
$\varphi:J_i\to J_j$ satisfies $u^i\varphi(m)=u^im\,\varphi(1)$, so
it is multiplication by $z=\varphi(1)=x+yv\in J_j$, with $x\in A$ and
$y\in u^{-j}A$. As $J_i$ is generated by $1$ and $v/u^i$, it remains
to impose
\[
 z\,\frac{v}{u^i}=y\,u^{\,n+1-i}+\frac{x}{u^i}\,v\in J_j .
\]
The coefficient of $v$ lies in $u^{-j}A$ exactly when
$x\in u^{i-j}A$, and the constant term lies in
$u^{\,n+1-i-j}A\subset A$ because $i+j\le2\delta\le n+1$. Hence
\[
 \Hom_R(J_i,J_j)=u^{\max(0,i-j)}A\oplus u^{-j}Av,
\]
and comparing with $S$ summand by summand gives
\begin{equation}\label{eq:an-hom-colength}
 \ell\bigl(S/\Hom_R(J_i,J_j)\bigr)
 =\max(0,i-j)+\delta-j
 =\delta+i-j-\min(i,j).
\end{equation}

\begin{corollary}[Sum-of-squares penalty]
\label{cor:an-explicit-penalty}
For the balanced type above, put $B_j=\sum_{i=j}^{\delta}b_i$.
Then
\begin{equation}\label{eq:an-explicit-penalty}
 c_R(M)
 =\sum_{i,j=1}^{\delta}b_i b_j\min(i,j)
 =\sum_{j=1}^{\delta}B_j^2.
\end{equation}
Consequently, on an integral projective complex curve of
arithmetic genus $g\ge2$ with one $A_n$ singularity, a general
rank-$r$, degree-$d$ bundle $E$ satisfies, for every saturated
sequence $0\to F\to E\to Q\to0$ of fixed kernel rank $k$,
\[
 D_e\ge\sum_{j=1}^{\delta}B_j^2,
 \qquad e=\deg F,
\]
where the $b_i$ are the multiplicities of the completed quotient
$\widehat Q$.
\end{corollary}
\begin{proof}
Replace $M$ by $\bigoplus_iJ_i^{\oplus b_i}$, whose normalization
is $S^a$. Summing~\eqref{eq:an-hom-colength} over the blocks of
its endomorphism matrix gives
\[
 \begin{aligned}
 \ell\bigl(\End_S(\bar M)/\End_R(M)\bigr)
 &=\sum_{i,j=0}^{\delta}b_i b_j
       \bigl(\delta+i-j-\min(i,j)\bigr)\\
 &=\delta a^2-
       \sum_{i,j=1}^{\delta}b_i b_j\min(i,j).
 \end{aligned}
\]
The $i-j$ terms cancel. Definition~\eqref{eq:ade-local-defect}
gives the first equality in~\eqref{eq:an-explicit-penalty};
the second follows from
\[
 \min(i,j)=\sum_{m=1}^{\delta}
          \mathbf1_{m\le i}\mathbf1_{m\le j}.
\]
The global bound is
Proposition~\ref{prop:ade-generic-geometry}
(\ref{part:ade-lower-bound}) applied to the quotient type.
\end{proof}

\begin{remark}
At the critical degree, $D_e=0$ forces
$b_1=\cdots=b_\delta=0$, recovering the boundary exclusion used
in the ADE counting theorem. At positive expected dimension,
the formula gives explicit necessary restrictions on occurrence.
For example, a quotient of completed type
$R^{s-b}\oplus I_i^{\oplus b}$, where $s=r-k$,
$1\le i\le\delta$ and $1\le b\le s$, requires
\[
 D_e\ge i b^2.
\]
For a node, $\delta=1$ and the penalty is $b_1^2$, which is the
codimension bound in \cite[Proposition~2.7]{Bho06}; in this case the bound
above is \cite[Proposition~2.8(4)]{Bho06}.
\end{remark}

\bibliographystyle{alpha}
\bibliography{ref}
\textsc{Department of Mathematics, University of Washington, Seattle, WA
98195, USA}

\textit{Email address: }\texttt{arka1996@uw.edu}
\end{document}